\documentclass{amsart}

\usepackage[small]{diagrams}
\usepackage{amssymb,amsmath,amsthm,xspace}

\def\mat#1{\ensuremath{#1}\xspace}

\def\makemath#1#2#3{
\newsavebox{#2}
\sbox{#2}{\ensuremath{#3}}
\def#1{\usebox{#2}\xspace}}

\def\cA{\mat{\mathbb{A}}}   

\def\cL{\mat{\mathbb{L}}}   %
\def\cQ{\mat{\mathbb{Q}}}   
\def\cR{\mat{\mathbb{R}}}   
\def\cC{\mat{\mathbb{C}}}   

\def\cZ{\mat{\mathbb{Z}}}   

\def\lE{\mat{\mathcal{E}}}

\def\lH{\mat{\mathcal{H}}}

\makemath{\Psi}{\Psibox}{\Psi}
\makemath{\Phi}{\Phibox}{\Phi}

\def\la{\mat{\lambda}}

\def\si{\mat{\sigma}}

\def\al{\mat{\alpha}}

\def\ksi{\mat{\xi}}

\def\mrm@#1{\mat{\mathrm{#1}}}

\def\DMO{\DeclareMathOperator}

\DMO{\Hom}{Hom}
\DMO{\lHom}{\lH\mathit{om}}
\DMO{\Ext}{Ext}
\DMO{\lExt}{\lE\mathit{xt}}
\DMO{\End}{End}
\DMO{\Aut}{Aut}
\DMO{\Fun}{Fun}
\DMO{\Tor}{Tor}
\DMO{\ext}{ext}

\DMO{\Ob}{Ob}
\DMO{\Mor}{Mor}
\DMO{\im}{im}
\DMO{\coim}{coim}
\DMO{\coker}{coker}
\DMO{\Arr}{Arr}

\DMO{\Id}{Id}
\DMO{\add}{add} 
\DMO{\ind}{ind} 
\DMO{\pro}{pro} 
\DMO{\Map}{Map}
\DMO{\Iso}{Iso}
\DMO{\Isom}{Isom}

\DMO{\Presh}{Presh}

\DMO\coalg{Coalg}
\DMO{\Rep}{Rep}

\DMO{\Mod}{Mod}
\DMO{\rad}{rad}
\DMO{\soc}{soc}
\DMO{\ann}{ann}

\DMO{\Spec}{Spec}
\DMO{\spec}{Spec}
\DMO{\Proj}{Proj}
\DMO{\supp}{supp}
\DMO{\Coh}{Coh}
\DMO{\coh}{Coh}
\DMO{\Qcoh}{QCoh}
\DMO{\QCoh}{QCoh}
\DMO{\Pic}{Pic}
\DMO{\Div}{Div}
\DMO{\ch}{ch}
\DMO{\Hilb}{Hilb}
\DMO{\Fitt}{Fitt}
\DMO{\Quot}{Quot}

\DMO{\Gras}{Gr}
\DMO{\Flag}{Flag}

\DMO{\cone}{cone}
\DMO{\Tw}{Tw}
\DMO{\rank}{rk}
\DMO{\rk}{rk}
\DMO{\codim}{codim}
\DMO{\cov}{cov}
\DMO{\sgn}{sgn}
\DMO{\td}{td}
\DMO{\GL}{GL}
\DMO{\SL}{SL}

\DMO\Der{Der}
\DMO\der{Der}
\DMO\coder{Coder}
\DMO{\diag}{diag}
\DMO{\HMod}{HMod} 
\DMO{\ad}{ad}
\DMO*{\colim}{colim}
\DMO*{\hocolim}{hocolim}
\DMO*{\holim}{holim}
\DMO{\Ho}{Ho}

\DMO{\har}{char}
\DMO{\sk}{sk}
\DMO{\cosk}{cosk}
\DMO{\Gal}{Gal}
\DMO{\tr}{tr}
\DMO{\Tr}{Tr}
\DMO{\Sh}{Sh}
\DMO{\Is}{Is} 
\DMO{\Hol}{Hol} 
\DMO{\Lie}{Lie} 
\DMO{\Res}{Res} 
\DMO{\irr}{irr} %
\DMO{\Irr}{Irr} %
\DMO{\Exp}{Exp} %
\DMO{\Log}{Log} %
\DMO{\mult}{mult} %
\DMO{\height}{ht}

\def\ts{\otimes}

\def\sb{\subset}

\def\xx{\times}

\def\lpoly#1{[\![#1]\!]} 

\def\arrowsUsual{
\def\arr{\xrightarrow}
\def\ar{\rightarrow}
\def\mto{\mapsto}
\def\emb{\hookrightarrow}    }

\newif\ifukr\ukrfalse
\newif\ifrus\rusfalse

\newtheorem{prp}{\ifukr Пропозиція \else \ifrus Предложение \else Proposition\fi\fi}[section]
\newtheorem{thr}[prp]{\ifukr Теорема \else \ifrus Теорема \else Theorem\fi\fi}

\theoremstyle{definition}

\newtheorem{exm}[prp]{\ifukr Приклад \else \ifrus Пример \else Example\fi\fi}
\newtheorem{rmr}[prp]{\ifukr Зауваження \else \ifrus Замечание \else Remark\fi\fi}

\arrowsUsual
\def\st{{}_{\rm st}}
\DMO{\wt}{wt}

\begin{document}

\title[]{Stringy motives of symmetric products}%

\author{Sergey Mozgovoy}%

\address{Institut f\"ur Mathematik, Johannes Gutenberg-Universit\"at Mainz,
55099 Mainz, Germany.}%

\email{mozgov@mathematik.uni-mainz.de}%

\thanks{}%
\subjclass[2000]{14Q15}%
\keywords{}%

\begin{abstract}
Given a complex smooth algebraic variety $X$, we compute the
generating function of the stringy motives of its symmetric powers
as a function of motive of $X$. In dimension two we recover the
G\"ottsche formulas for Hilbert schemes. We use the formalism of
$\lambda$-rings to get a particularly compact formula, which is
convenient for explicit computations.
\end{abstract}
\maketitle

\section*{Introduction}
The stringy motive $[X]\st$ of a log terminal complex algebraic
variety $X$ (see e.g.\ \cite{Bat1, Craw1}) is an invariant of $X$
such that whenever $X$ possesses a crepant resolution $Y\ar X$,
the motive $[Y]$ coincides with $[X]\st$. Here motives are
elements of certain extension $\hat A_\infty$ of the Grothendieck
ring $A$ of the category of Chow motives and the map $X{\mto}[X]$
sending an algebraic variety to its motive is the one constructed
by Gillet and Soul{\'e} \cite{GilSoule1}. The purpose of this
paper is to compute the stringy motives of symmetric products of
smooth manifolds.

Let us denote by $\cL$ the motive $[\cA^1]$, called the Lefschetz
motive. It is known that the Grothendieck ring $A$ has a structure
of the \la-ring, where the \si-operations are defined on the level
of varieties by symmetric products (see e.g.\ \cite{Getz1} for
this result and basic definitions concerning \la-rings). In
particular, $\si_n(\cL)=\cL^n$. Given any \la-ring $R$, we endow
the ring $R\lpoly t$ with a structure of \la-ring by
$\si_n(at^k):=\si_n(a)t^{nk}$, $a\in R$. Define the map
$\Exp:R\lpoly t^+\ar 1+R\lpoly t^+$ (see e.g. \cite{Getz1} or
\cite[Appendix]{M2}) by
$$\Exp(f):=1+\si_1(f)+\si_2(f)+\dots$$
Clearly $\Exp(f+g)=\Exp(f)\Exp(g)$. Using the Adams operations
(and assuming that $\cQ\sb R$) one can also write
$$\Exp(f)=\exp\Big(\sum_{k\ge1}\frac 1k\psi_k(f)\Big).$$
Our main result is the following

{\def\theprp{1}
\begin{thr}\label{thr intro main}
Let $X$ be a complex smooth algebraic variety of dimension $d\ge
2$. Then
$$\sum_{n\ge 0}{[X^{(n)}]\st}t^n
=\Exp\left(\frac {[X]\, t}{1-\cL^{d/2}t}\right),$$ %
where $X^{(n)}$ denotes the symmetric product.
\end{thr}}

If $X$ is a surface then there is a canonical crepant resolution
$X^{[n]}\ar X^{(n)}$, where $X^{[n]}$ denotes the Hilbert scheme
of points and we get
$$\sum_{n\ge 0}{[X^{[n]}]}t^n
=\Exp\left(\frac {[X]\, t}{1-\cL t}\right).$$

This formula is equivalent to the formula of G\"ottsche
\cite{Go1,Go2,GS1,Cheah1,Go3} but is more compact and probably
more convenient for explicit computations (certainly, using some
computer algebra system). In such computations one can use the
above formula for $\Exp$ in terms of the Adams operations as the
Adams operations on the Hodge (or Poincar\'e) polynomials are
particularly easy to describe. We discuss this together with some
examples in the second section.

The idea of the proof of the theorem is to use the result of
Batyrev \cite{Bat1} saying that the stringy motive of a finite
quotient of a smooth manifold coincides with an orbifold motive.
In our situation, we consider the action of the symmetric group
$S_n$ on $X^n$, where the orbifold motive is very easy to compute.
It should be noted that originally the Euler number of the Hilbert
scheme was computed first and then it was shown that it coincides
with the orbifold Euler number \cite{HH1}. Thus, our proof is
rather artificial but is probably the shortest one.

\section{Stringy motives}
For the definitions and basic properties of stringy motives and
motivic integration we refer to \cite{Bat1,Craw1}. We will only
define the ring where our motives live. Let $A$ be the
Grothendieck ring of Chow motives. It is an algebra over
$\cZ[\cL]$ (here \cL is considered as a variable). We define
$A_n:=A\ts_{\cZ[\cL]}\cZ[\cL^{\pm 1/n}]$ and $A_\infty:=\colim
A_n$; thus, $A_\infty$ contains arbitrary rational powers of \cL.
There is a non-Archimedean norm $|-|:A_\infty\ar\cR_{\ge0}$
satisfying $|[X]|=2^{\dim X}$ and $|\cL^s|=2^s$. Let $\hat
A_\infty$ be the completion of $A_\infty$ with respect to this
norm. As it was mentioned in the introduction, the ring $A$ has
the structure of a \la-ring; it can be extended to $\hat A_\infty$
in such a way that $\si_n(a\cL^s)=\si_n(a)\cL^{ns}$ for any $a\in
A$ and $s\in\cQ$. For any log terminal algebraic variety $X$, one
can define an invariant $[X]\st\in\hat A_\infty$ called the
stringy motive of $X$ (see \cite{Bat1,Craw1}).


Let $X$ be a complex smooth algebraic variety of dimension $d$,
$G$ be a finite group acting effectively on $X$ and suppose that
there exists a geometric quotient $X/G$. We assume that the
ramification locus, that is $\cup_{g\ne1}X^g$, has codimension at
least $2$ in $X$ and that $X^g$ is connected for any $g\in G$. Let
us define the orbifold motive of the pair $(X,G)$ according to
\cite{Bat1}. Given a vector space $V$ of dimension $d$ and a
linear operator $T$ on $V$ of finite order, let $(e^{2\pi
i\al_1},\dots,e^{2\pi i\al_d})$ with $0\le\al_k<1$ be the
eigenvalues of $T$. Define the weight
$$\wt(T):=\sum_{k=1}^d\al_k\in\cQ.$$
For any $g\in G$, define the weight $\wt(g)$ to be the weight of
the action of $g$ on $T_xX$, for some $x\in X^g$. It does not
depend on the choice of the point $x\in X^g$. Define then the
orbifold motive
$$[X/G]_{\rm orb}:=\sum_{[g]\sb G}[X^g/C(g)]\cdot\cL^{\wt(g)},$$
where $[g]$ runs over the conjugacy classes of $G$ (or rather
their representatives) and $C(g)\sb G$ is the centralizer of $g$.

\begin{thr}[{see \cite[Theorem 7.5]{Bat1}}]
It holds $[X/G]\st=[X/G]_{\rm orb}$.
\end{thr}

We apply this theorem in the case of the symmetric group $S_n$
acting in a natural way on $X^n$, where $X$ is smooth. For any
partition $\la=(\la_1,\la_2,\dots)=(1^{a_1}\dots r^{a_r})$, we
define $X^{\la}:=X^{(a_1)}\xx\dots\xx X^{(a_r)}$. We define the
length of \la to be $l(\la):=a_1+\dots+a_r$ and the weight of \la
to be $|\la|:=\sum_{i\ge 1}\la_i=\sum_{k=1}^r ka_k$.

\begin{prp}
Let $X$ be a complex smooth algebraic variety of dimension
$d\ge2$. Then, for any $n\ge 1$,
$$[X^{(n)}]\st=\sum_{|\la|=n}[X^{(\la)}]\cdot \cL^{\frac d2(n-l(\la))},$$
where the sum runs over all partitions \la of $n$.
\end{prp}
\begin{proof}
Consider the natural action of $S_n$ on $Y=X^n$. The conjugacy
classes of $S^n$ are parameterized by the partitions of $n$, where
a partition $\la=(1^{a_1}\dots r^{a_r})$ corresponds to the
permutations $g\in S_n$ that have $a_1$ $1$-cycles, $a_2$
$2$-cycles etc. The manifold $Y^g$ is isomorphic to
$X^{a_1}\xx\dots X^{a_r}$ and has codimension at least $d\ge 2$ in
$Y$ whenever $\la\ne(1^n)$, i.e., $g\ne 1$. The cenralizer of $g$
consists of those elements of $S_n$ that permute the cycles of
$g$. It follows that the variety $X^g/C(g)$ is isomorphic to
$X^{(a_1)}\xx\dots\xx X^{(a_n)}=X^{(\la)}$. We note that the cycle
$(12\dots k)$ acting by permutation on $\cC^k$ has an eigenvalue
$\ksi=e^{2\pi ij/k}$ with an eigenvector
$(\ksi,\ksi^2,\dots,\ksi^k)$ for any $j=0,\dots, k-1$. This
implies that the weight of the permutation equals
$(1+\dots+(k-1))/k=(k-1)/2$. It follows that the cycle $(12\dots
k)$ acting on $(\cC^d)^k$ by permutation has weight $d(k-1)/2$.
Finally, an element $g$ of type $\la=(1^{a_1}\dots r^{a_r})$ in
$S_n$ has weight
$$\wt(g)=\sum_{k=1}^r \frac {d(k-1)}2 a_k=\frac d2(n-l(\la)).$$
Now the statement of the proposition follows from the previous
theorem.
\end{proof}

We are ready to prove the main result of the paper.


\begin{proof}[Proof of Theorem \ref{thr intro main}]
Using the previous proposition we get
\begin{align*}
\sum_{n\ge 0}[X^{(n)}]\st t^n%
  =&\sum_\la [X^{(\la)}]\cdot \cL^{\frac d2(|\la|-l(\la))}\cdot t^{|\la|}\\
  =&\sum_{\la=(1^{a_1}\dots r^{a_r})}
    \prod_{k=1}^r([X^{(a_k)}]\cdot\cL^{\frac d2 (ka_k-a_k)}\cdot t^{ka_k})\\
  =&\sum_{\la=(1^{a_1}\dots r^{a_r})}
    \prod_{k=1}^r\si_{a_k}([X]\cdot\cL^{\frac d2(k-1)}\cdot t^k)\\
  =&\prod_{k\ge 1}\sum_{a\ge0}\si_{a}([X]\cdot\cL^{\frac d2(k-1)}\cdot t^k)
  =\prod_{k\ge 1}\Exp([X]\cdot\cL^{\frac d2(k-1)}\cdot t^k)\\
  =&\Exp(\sum_{k\ge 0}[X]\, t(\cL^{d/2}t)^k)
  =\Exp\left(\frac {[X]\, t}{1-\cL^{d/2}t}\right).
\end{align*}
\end{proof}

\begin{rmr}
In the same way as above one can show also
$$\sum_{n\ge 0}\frac{[X^{(n)}]\st}{\cL^{dn/2}}t^n
=\Exp\left(\frac {[X]}{\cL^{d/2}}\frac t{1-t}\right).$$ %
\end{rmr}

\section{Examples}
We endow the ring $\cZ[x_1,\dots, x_r]$ with the structure of a
\la-ring by $\si_n(x^\al)=x^{n\al}$, where $\al\in\cZ_+^r$. Most
conveniently, this \la-structure can be written using the Adams
operations
$$\psi_n(f(x_1,\dots,x_r))=f(x_1^n,\dots,x_r^n),\qquad f\in\cZ[x_1,\dots,x_n].$$
In order to avoid any problems with the Adams operations in what
follows, we tensor our \la-rings with \cQ without mentioning that
and so we always assume that our \la-rings contain \cQ.

There are three basic realizations of the ring of motives $A$. The
Euler number $e:A\ar\cZ$, the virtual Poincar\'e polynomial
$P:A\ar \cZ[v]$, and the virtual Hodge polynomial (also called
Hodge-Deligne polynomial) $E:A\ar\cZ[u,v]$ (see e.g.\
\cite{DanKhov}). The Euler number of a variety can be defined as
the Euler characteristic of the complex of cohomologies with
compact support. To define the virtual Hodge polynomial or the
virtual Poincar\'e polynomial, one uses the mixed Hodge structure
on the cohomologies with compact support. In the case of the
Lefschetz motive one has $e(\cL)=1$, $P(\cL)=v^2$, $E(\cL)=uv$.
The relations between the above three realizations are
$$e(X)=P(X; 1),\qquad P(X; v)=E(X; v,v).$$
An important property of the above realizations is that they are
morphisms of \la-rings. In the case of the Poincar\'e polynomial
this corresponds to the results of Macdonald \cite{Mac2}. In the
case of the Hodge polynomial see e.g.\ \cite{Ste1}. All three
realizations can be extended to $\hat A_\infty$ after extending
accordingly the rings $\cZ$, $\cZ[v]$, and $\cZ[u,v]$. In this way
one can also define the stringy Euler number, the stringy
Poincar\'e function, and the stringy Hodge function (usually
called stringy E-function). After all these remarks we can rewrite
Theorem \ref{thr intro main} substituting there any of three
described realizations. For example, given a complex smooth
algebraic variety $X$ of dimension $d\ge 2$, it holds
$$\sum_{n\ge 0}{P(X^{(n)})\st}t^n
=\Exp\left(\frac {P(X)\, t}{1-v^dt}\right).$$ %

\begin{exm}
For the Euler numbers we get
$$\sum_{n\ge 0}{e(X^{(n)})\st}t^n
=\Exp\left(\frac {e(X)\,t}{1-t}\right)%
=\prod_{k\ge1}\Exp(t^k)^{e(X)}%
=\prod_{k\ge1}(1-t^k)^{-e(X)}.%
$$ %
For surfaces, this is the G\"ottsche formula for the Euler numbers
of Hilbert schemes.
\end{exm}

As it has already been mentioned, the \la-structure on the ring of
polynomials is written most conveniently using the Adams
operations. In view of the above formulas, it is important to know
the formula for $\Exp$ in terms of the Adams operations. There is
a formal identity between symmetric functions (see e.g.
\cite[2.10]{Mac1})
$$\sum_{n\ge0}h_nt^n=\exp(\sum_{k\ge1}\frac 1kp_kt^k)$$
which implies
$$\Exp(f)=\sum_{n\ge 0}\si_n(f)=\exp(\sum_{k\ge1}\frac 1k\psi_k(f)).$$

\begin{exm}
Let us find $P(X^{[2]})$, where $X$ is a K3-surface. It is known
that $P(X;v)=v^4+22v^2+1$. Applying the main theorem to the
stringy Poincar\'e functions we get, noticing that the Adams
operations are ring homomorphisms
\begin{align*}
&1+{P(X)}t+{P(X^{[2]})}t^2+o(t^2)%
=\Exp\left({P(X)}t(1+v^2t)\right)+o(t^2)\\
=&\exp\big({P(X)}(t+v^2t^2)+\frac12{\psi_2(P(X))}t^2\big)+o(t^2).\\
\end{align*}
This implies
$$P(X^{[2]})=v^2 P(X)+\frac 12\psi_2(P(X))+\frac 12P(X)^2=
v^8+23v^6+276v^4+23v^2+1.$$
\end{exm}



\bibliography{fullbib}
\bibliographystyle{hamsplain}

\end{document}